\documentclass[12pt]{article}        
\usepackage[french]{babel}     

\newtheorem{thm}{Th\'eor\`eme}
\newtheorem{prop}{Proposition}

{\itshape}{\rm}

\let\bs=\bigskip

\def\R{{\msb R}}
\def\C{{\msb C}}

\def\F*2g{{\msb F}^*_{2^g}}
\def\f#1{{\msb F}_{#1}}

\font\tenmsb=msbm10 at 12pt
\font\sevenmsb=msbm9
\font\fivemsb=msbm6
\newfam\msbfam%
\def\msb{\fam\msbfam\tenmsb}% 
\textfont\msbfam=\tenmsb \scriptfont\msbfam=\sevenmsb%
\scriptscriptfont\msbfam=\fivemsb%

\def\val{\mathop {\rm val}}

\begin{document}

\setcounter{prop}{12}

\title{Errata \`a ``Sur les repr\'esentations non ramifi\'ees des groupes
r\'eductifs
$p$-adiques;\\ l'exemple de ${\rm GSp}(4)$''}
\author{Fran\c cois Rodier\\Institut de Math\'ematiques de Luminy --
C.N.R.S.\\Marseille -- France}
\date{}
\maketitle

\begin{abstract}\noindent
On corrige deux erreurs dans le papier \cite{ro}: l'une dans l'\'etude d'une
involution sur les repr\'esentations irr\'eductibles non ramifi\'ees d'un groupe
semi-simple, l'autre dans la description de repr\'esentations du groupe
${\rm GSp}(4)$.

\bs%\bs
\centerline{\bf Abstract}\noindent
We correct two errors in the paper \cite{ro}: the first in the study of an
involution on the irreducible unramified representations of a semi-simple
group, the second in the description of representations of the group ${\rm
GSp}(4)$.

\end{abstract}

\bs
Deux erreurs m'ont \'et\'e signal\'ees dans cet article, la premi\`ere par
Amritanshu Prasad, qui avait utilis\'e l'\'enonc\'e de la proposition 13  dans son
papier
\cite{ap1} et qui 
a d\^u \'ecrire par la suite  un erratum, et la seconde par Laurent Clozel. Je
les remercie tous deux de m'avoir signal\'e ces erreurs.

\section{L'erreur dans la Proposition 13}
%====================================

\def\sgn{\mathop {\rm sgn}}
\def\Int{\mathop {\rm Int}}

La premi\`ere erreur concerne la proposition 13.
Elle a des cons\'equences sur la description des composants
irr\'eductibles des repr\'esentations de ${\rm GSp}(4)$ dans le chapitre 6 et
dans la remarque finale du chapitre 7, mais ni sur le nombre de ces
repr\'esentations, ni sur leur multiplicit\'e. 

Elle est due \`a une confusion entre deux notations.
La notation $\sgn$ d\'efinie dans la section  2.1 comme d\'enotant un
carac\-t\`ere de
$k^\times$, donc une application du groupe multiplicatif du corps local non
archim\'edien
$k$ dans $\C$, est \`a ne pas confondre avec la notation
$\sgn q_x$ d\'efinie en 5.3, o\`u $q_x$ repr\'esente un volume, donc un nombre r\'eel.

Dans la d\'emonstration de la proposition 13, l'assertion
$q_t=\rho_P(t)^{-2}$ \'etait utilis\'ee pour prouver que $\sgn q_t=1$, alors que
$\rho(t)$ n'est pas forc\'ement entier puissance de $q$.

L'assertion de la ligne suivante, obtenue \`a l'aide du lemme 4, doit s'\'ecrire par
cons\'equent:
$$((\hat\pi)_U(t)\circ A)x_U=\sgn q_t\
A(\pi((T_{w_0(t^{-1})})^{-1})x)_U.$$

Elle implique 
$$(\hat\pi)_U(t)\circ A=\sgn q_t\ \rho_P(t)^{2}
A\circ\pi_U(w_0(t))$$

Ou encore, en remarquant que 
$\sgn q_t=sgn\ \rho_P^2(t)=\sgn\nu(t)$:
$$(\hat\pi)_U(t)\circ A=\sgn\circ\nu(t) \rho_P(t)^{2}
A\circ\pi_U(w_0(t))$$

D'ou l'\'enonc\'e corrig\'e de la proposition 13:

\begin{prop}
. --- La repr\'esentation $R(\hat \pi)$ est \'equivalente \`a
$(\sgn\circ\nu) R(\pi)\circ \Int w_0$ o\`u $\Int w_0$ est l'automorphisme de $T$
d\'efini par
$w_0$.
\end{prop}

Par la suite cette erreur  affecte  dans les sections post\'erieures la
description des composants irr\'eductibles  du groupe ${\rm GSp}(4)$. 

Dans la section 6.2, il faut lire que la repr\'esentation $(\sgn\circ\nu)
R(\hat\pi_\chi)$ est compos\'ee de 
$\chi$, avec la multiplicit\'e 2, et de $w_\alpha \chi$ avec la multiplicit\'e 1.
Et par cons\'equent $I(\chi)$ est compos\'ee de 4 repr\'esentations irr\'eductibles:
$\pi_\chi$,
$\pi'_\chi$, $(\sgn\circ\nu) \hat\pi_\chi$, $(\sgn\circ\nu) \hat\pi'_\chi$.

Dans la section 6.3:  $I(\chi)$ est compos\'ee de 2 repr\'esentations
irr\'eductibles: 
$\pi_\chi$,
$(\sgn\circ\nu) \hat\pi_\chi$.

C'est aussi le m\^eme cas dans les sections 6.4 et 6.5.

A la fin du paragraphe 7.2, il faut modifier la remarque finale:

{\it Remarque.} --- Si $\chi$ est le caract\`ere
$t\mapsto|t^{2\alpha+\beta}|^{1/2}$ ou le caract\`ere
$t\mapsto\sgn\circ\nu(t)|t^{2\alpha+\beta}|^{1/2}$, alors
$I(\chi)$ admet $\sgn\circ\nu\ \hat\pi_\chi$ et $\pi'_\chi$
comme composants temp\'er\'es.

\section{Erreurs dans l'\'enonc\'e du th\'eor\`eme 2}

Les conditions sur le repr\'esentant $x$ dans $X(T)\otimes \R$ du caract\`ere
unitaire
$\chi_U$ impos\'ees dans la section 7.2
sont  traduites maladroitement dans le th\'eor\`eme 2. Voici les corrections.

 \setcounter{thm}{1}

\begin{thm}
. --- Les repr\'esentations irr\'eductibles non ramifi\'ees de $G$
sont les suivantes.

\begin{description}
\item[\hbox{\rm a)}] sans changement;
\item[\hbox{\rm b)}] 
$I(\chi)$ pour $\chi(t)=\exp(-2\pi i\mu\val t^{\alpha+\beta})
|t^{\alpha}|^\lambda$ avec $\lambda,\mu\in\R$, $0<\mu<1/2$ et 
 $0<\lambda<1/2$;
\item[\hbox{\rm c)}]
sans changement;
\item[\hbox{\rm d)}]
$I(\chi)$ pour $\chi(t)=\exp(-2\pi i\mu\val t^{2\alpha+\beta})
|t^{\beta}|^\lambda$ avec $\lambda,\mu\in\R$, $0<\mu<1/2$ et 
 $0<\lambda<1/2$;
\item[\hbox{\rm e)}]
sans changement;
\item[\hbox{\rm f)}]
$I(\chi)$ pour $\chi(t)=\sgn t^{\alpha+\beta}\exp(-2\pi i\mu\val t^{\alpha})
|t^{\alpha+\beta}|^\lambda$ avec $\lambda,\mu\in\R$,
\mbox{$0<\mu<1/2$} et 
 $0<\lambda<1/2$;
\item[\hbox{\rm g)}]
sans changement;
\item[\hbox{\rm h \`a n)}]
sans changement;
\item[\hbox{\rm o)}]
Les composants de $I(\chi)$ pour 
$\chi=|t^\alpha|^\lambda |t^\beta|^{1/2}$ avec $1/2\le \lambda\le 1$;
\item[\hbox{\rm p)}]
Les composants de $I(\chi)$ pour 
$\chi=\sgn\circ\nu|t^\alpha|^\lambda |t^\beta|^{1/2}$ avec $1/2\le
\lambda\le 1$.
\end{description}

\end{thm}

%+++++++++++++++++++++++++++++++++++

%+++++++++++++++++++++++++++++++++++
\end{document}